\documentclass[a4paper, 10pt, twocolumn]{article}

\usepackage{amsfonts,amssymb,amsmath,amsthm}
\usepackage{subfigure}
\usepackage{graphicx}
\usepackage[footnotesize]{caption}

\renewcommand{\thepage}{}
\usepackage[top=1.5cm, left=1.5cm, right=1.5cm, bottom=1.5cm]{geometry}

\title{SLS (Single $\ell_1$ Selection): a new greedy algorithm with an $\ell_1$-norm selection rule}

\author{Ramzi Ben Mhenni$^1$, S\'ebastien Bourguignon $^1$ and J\'er\^ome Idier$^{1}$.\\
\footnotesize $^1$ \'Ecole Centrale de Nantes - CNRS, Laboratoire des Sciences du Num\'erique de Nantes, 1 rue de la No\"e, F-44321 Nantes, France.}
\date{\empty} % no need for a date

\renewenvironment{abstract}{\bf\small {\em\ Abstract---}}{}

\newcommand{\vect}[1]{\boldsymbol{#1}}

\newcommand{\vx}{\vect{x}}

\newcommand{\vy}{\vect{y}}

\newcommand{\vr}{\vect{r}}

\newcommand{\vh}{\vect{h}}

\newcommand{\va}{\vect{a}}

\newcommand{\mat}[1]{\mathbf{#1}}
\newcommand{\mH}{\mat{H}}

\newcommand{\mA}{\mH_{\Qbar}}
\newcommand{\mB}{\mH_{\Qup}}

\newcommand{\mI}{\mat{I}}
\newcommand{\mP}{\mat{P}}

\newcommand{\Q}{\text{$\mathbb{S}$}}
\newcommand{\Qbar}{\overline{\Q}}

\newcommand{\Qup}{\Q_1}

\newcommand{\nup}{{n_1}}

\newcommand{\barn}{\bar{n}}

\newcommand{\n}{M}
\newcommand{\N}{N}
\newcommand{\Rn}{\R^\n}
\newcommand{\RN}{\R^\N}%{\mathbb{R}^N}
\newcommand{\R}{\mathbb{R}}
\newcommand{\norm}[1]{\|#1\|}

\def\sansprofil#1{}%#1

\def\argmin{\mathop{\mathrm{argmin}}}
\def\argmax{\mathop{\mathrm{argmax}}}
\renewcommand{\Q}{{\mathbb{S}}}
\renewcommand{\Qup}{\Q}
\renewcommand{\nup}{k}
\renewcommand{\barn}{\bar{k}}

\def\score{\mathcal{F}}
\usepackage{url}
\usepackage{color}
\usepackage[linesnumbered,ruled,vlined,english,onelanguage]{algorithm2e}

\usepackage{pgfplots,pgfplotstable}
\pgfplotsset{width=7cm,compat=1.8}
\usetikzlibrary{math}
\usetikzlibrary{positioning}
\usetikzlibrary{decorations.pathmorphing, shadows}
\usetikzlibrary{calc}

\def\mH{\mat{A}}
\def\vh{\vect{a}}
\begin{document}

\maketitle

\begin{abstract}
	In this paper, we propose a new greedy algorithm for sparse approximation, called SLS for {\it Single $\ell_1$ Selection}. 
	SLS essentially consists of a greedy forward strategy, where the selection rule of a new component at each iteration is based on solving a least-squares optimization problem,  penalized by the $\ell_1$ norm of the remaining variables. Then, the component with maximum amplitude is selected. Simulation results on difficult sparse deconvolution problems involving a highly correlated dictionary reveal the efficiency of the method, which outperforms popular greedy algorithms and Basis Pursuit Denoising  when the solution is sparse. 
\end{abstract}

\section{Introduction}
\label{sec:intro}

We consider the cardinality-constrained least-squares problem:
\begin{equation}
\min_{\vx} { \tfrac{1}{2} {\norm{\vy-\mH\vx}}^2_2 \;\text{ subject to (s.t.) } \; \norm{\vx}_0 \leq K },
\label{eq:P0}
\end{equation}
where $\vy \in \RN$, $\vx\in \Rn$, $\mH$ is a known matrix or {\it dictionary}, and 
 $\norm{\vx}_0$ is the $\ell_0$ ``norm'': $\norm{\vx}_0:= \text{Card} \{x_j| x_j \neq 0\}$. 
Problem~\eqref{eq:P0} is encountered in many sparse approximation problems occurring in inverse problems~\cite{Taylor79,Mendel83,Zala92}, denoising,  compression~\cite{Eldar12book}, or subset selection in Statistics~\cite{Miller02}.

Finding such best $K$-sparse solution  is NP hard~\cite{natarajan_sparse_1995}, therefore most works in signal processing and statistics have concentrated on developing computationally efficient, suboptimal, algorithms~\cite{Tropp10}. 
Forward greedy algorithms such as Orthogonal Matching Pursuit (OMP) \cite{Pati93} and Orthogonal Least Squares (OLS) \cite{Chen89} iteratively add new components to an initially empty model, then providing a $K$-sparse approximation in no more than $K$ iterations. However, the selection step of a new atom in such methods is highly sensitive to 
interferences between the dictionary atoms, in particular in the case of highly correlated dictionaries~\cite{soussen_bernoulli_2011}. 
Convex optimization strategies in which the $\ell_0$ norm in problem~\eqref{eq:P0} is replaced by the $\ell_1$ norm  $\norm{\vx}_1 := \sum_n |x_n|$ (which is known as the LASSO in Statistics~\cite{Tibshirani96} ), is another widespread approach, for which many dedicated algorithms have been proposed in the past twenty years. Optimizing all variables together in a convex approach then brings more robustness toward the aforementioned interferences, but the solution may often contain undesired nonzero components with small amplitudes. 

%In the compressed sensing framework~\cite{Eldar12book}, conditions on $\mH$ have been established, for which such approaches are ensured to solve the initial problem~\eqref{eq:P0}. However, in the case of inverse problems where $\mH$ is ill-conditioned, or if the mutual correlation of its columns is high, such guarantees are lost and, in practice, they often fail in finding the global optimum~\cite{Bourguignon2016}. 

%While the $\ell_1$-norm convex formulation exactly solves an approximate problem, the second class of methods returns a local optimum for the exact problem. The relative performance of methods is then evaluated under a compromise criterion between the quality of the solution and its computational cost. 

In this paper, we propose an algorithm which gathers advantages of the two classes of methods. It essentially consists of a greedy strategy, where the selection rule at each iteration is based on exploiting $\ell_1$-norm solutions. 
The number of iterations is then controlled by the sparsity level $K$ of the searched solution, limiting the computational burden. Moreover, the selection of each new atom, based on solving a convex optimization problem, is expected to be more robust to interferences between the different atoms than standard greedy methods. 

%More precisely, at a given iteration, an $\ell_1$-norm-penalized problem is built, where the $\ell_1$ norm operates on the non-selected variables, and the regularization  path is computed by an homotopy continuation algorithm~\cite{Tibshirani96,Donoho08}.A heuristic selection rule is then proposed, which considers the predominant variable along the regularization path. 

\section{Forward greedy algorithms}
\label{sec:glouton}

Forward greedy methods start from an empty set and iteratively construct a sparse solution. %, by alternating between two steps: a new atom is selected by maximizing a score function, denoted by $\score$, and then the model is updated. 
Let $\Qup$ denote the index set of the variables already selected (the current {\it support} of the solution, with $\nup$ components) and let $\Qbar$ index the remaining $\barn$ variables. 
In the following, $\mB$  denotes the matrix composed of the columns of $\mH$ indexed by $\Qup$. Similarly, $\vx_{\Qup}$ is the vector collecting the elements of $\vx$ indexed by $\Qup$. composed of the columns of $\mH$ (resp. the elements of $\vx$) indexed by $\Qup$. 
The principle of forward selection algorithms is given in  Algorithm~\ref{algo:glouton}.
\begin{algorithm}[h]			
	\SetAlgoLined
	\textbf{Initialization:} \quad $\Qup = \emptyset $\\% \quad \vr = \vy$ \\
	%\hrulefill \\
	\While{ $|\Qup|<K$}{
		%	~\\
		Variable selection: 
		$\hat{\jmath} = \argmax_{j \in \Qbar}~\score(j)$ ; \\%[0.2cm]
		Support update: $\Qup \leftarrow \Qup \cup \{\hat{\jmath}\}$; \\
		%	$\Qup \leftarrow \Qup \cup \{\hat{\jmath}\}$;\qquad \% \textcolor{blue}{ Support update} \% \\[0.2cm]
		(If needed) update the estimate $\vx_{\Qup}$ and the residual $ \vr = \vy- \mB \vx_{\Qup} $ ; \\%[0.2cm]
		%	$\vr = \vy- \mB  \vx_{\Qup} ;$\qquad \% \textcolor{blue}{ Residual update} \% 					
	}
	\KwResult{ support $\Qup$ and solution $\vx_{\Qup}$}
	\caption{Forward Selection greedy algorithm}
	\label{algo:glouton}
\end{algorithm}

In the sequel, we suppose that all columns in $\mH$ have unit norm. 
For OMP, the selected atom is the most correlated to the residual:
\begin{equation}
\label{eq:scoreOMP}
\score_{\text{OMP}} (j) = |\vh_j^T \vr|, \; j \in \Qbar,
\end{equation}
where $\vh_j$ is the $j$-th column of $\mH$. 
OMP includes an additional orthogonalization step of the solution on its support by:
$$
{\vx}_{\Qup} = \argmin_{\vx_{\Qup} \in \mathbb{R}^\nup} {\norm{\vy - \mB \vx_{\Qup}}}_2^2=
\mB^+ \vy,
$$
where $\mB^+ := {(\mB^T \mB)}^{-1} \mB^T$ denotes the pseudo-inverse of $\mB$.
For OLS \cite{Chen89}, the approximation error is minimized among all possible supports including one new component:
$$\widehat{\jmath} = \argmin_{j\in \Qup}\min_{\vx_{\Qup} \in \R^{\nup}} {\norm{\vy- \mH_{\Qup \cup \{j\}} \vx_{\Qup \cup \{j\}}}}_2^2,$$ 
which amounts to%, due to orthogonality, to
%$$\score_{\text{OLS}} (j) = \norm{\mH_{\Qup \cup \{j\}}(\mH_{\Qup \cup \{j\}}^T\mH_{\Qup \cup \{j\}})^{-1}\mH_{\Qup \cup \{j\}}^T \vy}^2.$$
$$\score_{\text{OLS}} (j) = \norm{\mH_{\Qup \cup \{j\}}\mH_{\Qup \cup \{j\}}^+ \vy}^2.$$

Restricting the selection step to models involving no more than one new component is a major limitation of such greedy algorithms. 
Let us consider the sparse deconvolution problem, where $\mH$ is composed of shifted versions of the impulse response of the filter. 
In the toy example of Figure\,\ref{fig::ex}, $\vx$ is composed of two close spikes, giving strongly overlapping echoes in the data $\vy$. The score function for the first iteration of both OMP and OLS is $\score(j) = |\vh_j^T \vy|$, and is shown in Figure\,\ref{fig::ex}\,(c). It is maximal for the index located in the middle of the two true indices, thus selecting a wrong atom. 
%and observe the selection of the first variable given by the OMP and OLS. See the example a) in Figure\,\ref{fig::ex} where

\def\wth{27mm}
\begin{figure}[h]
	\centering
	\footnotesize
	\begin{tabular}{@{}c@{\hspace*{0.15cm}}c@{\hspace*{0.15cm}}c@{}}
		(a) & (c) & (e)\\
		\includegraphics[width=\wth]{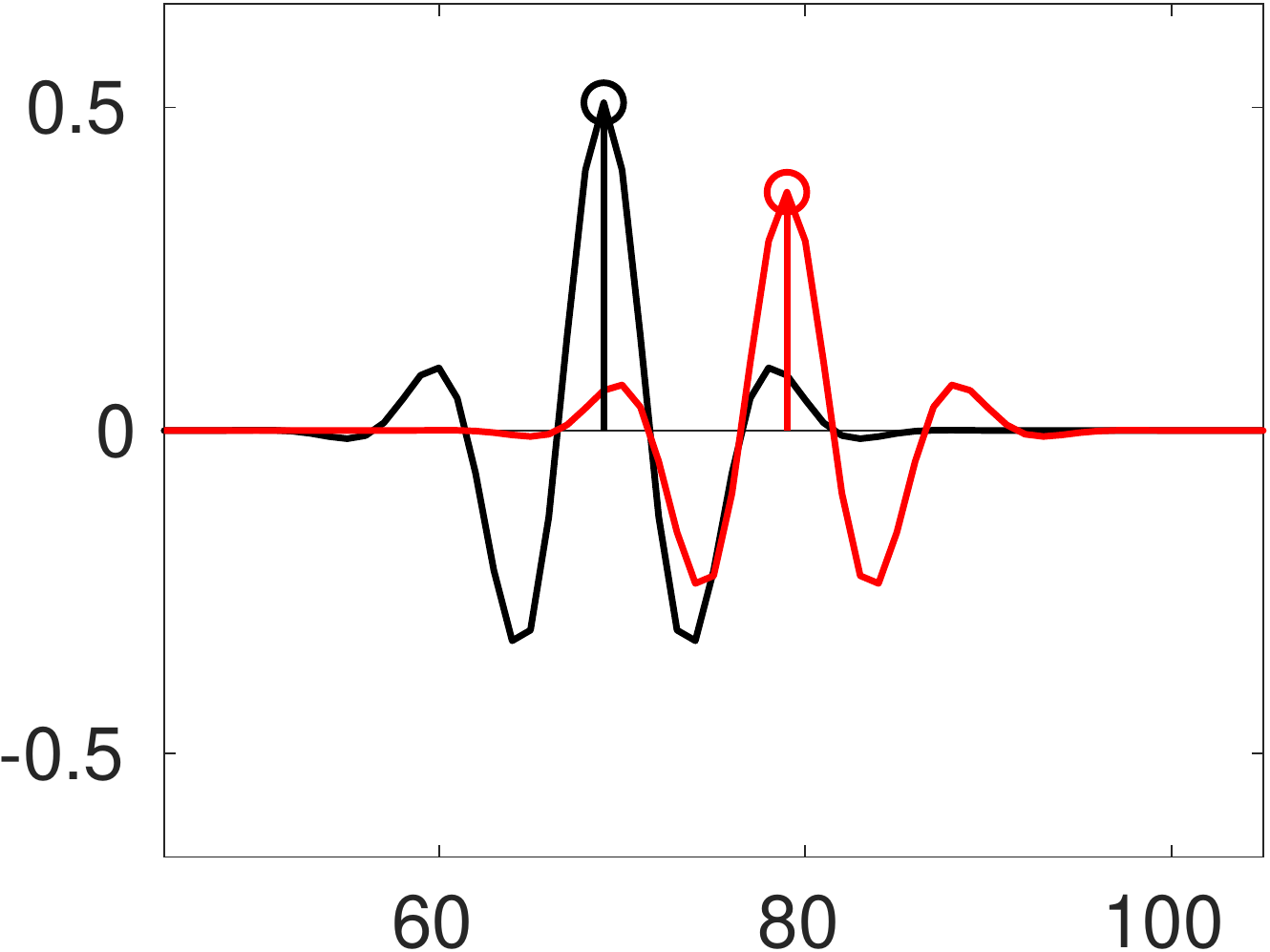}
		&	\includegraphics[width=\wth]{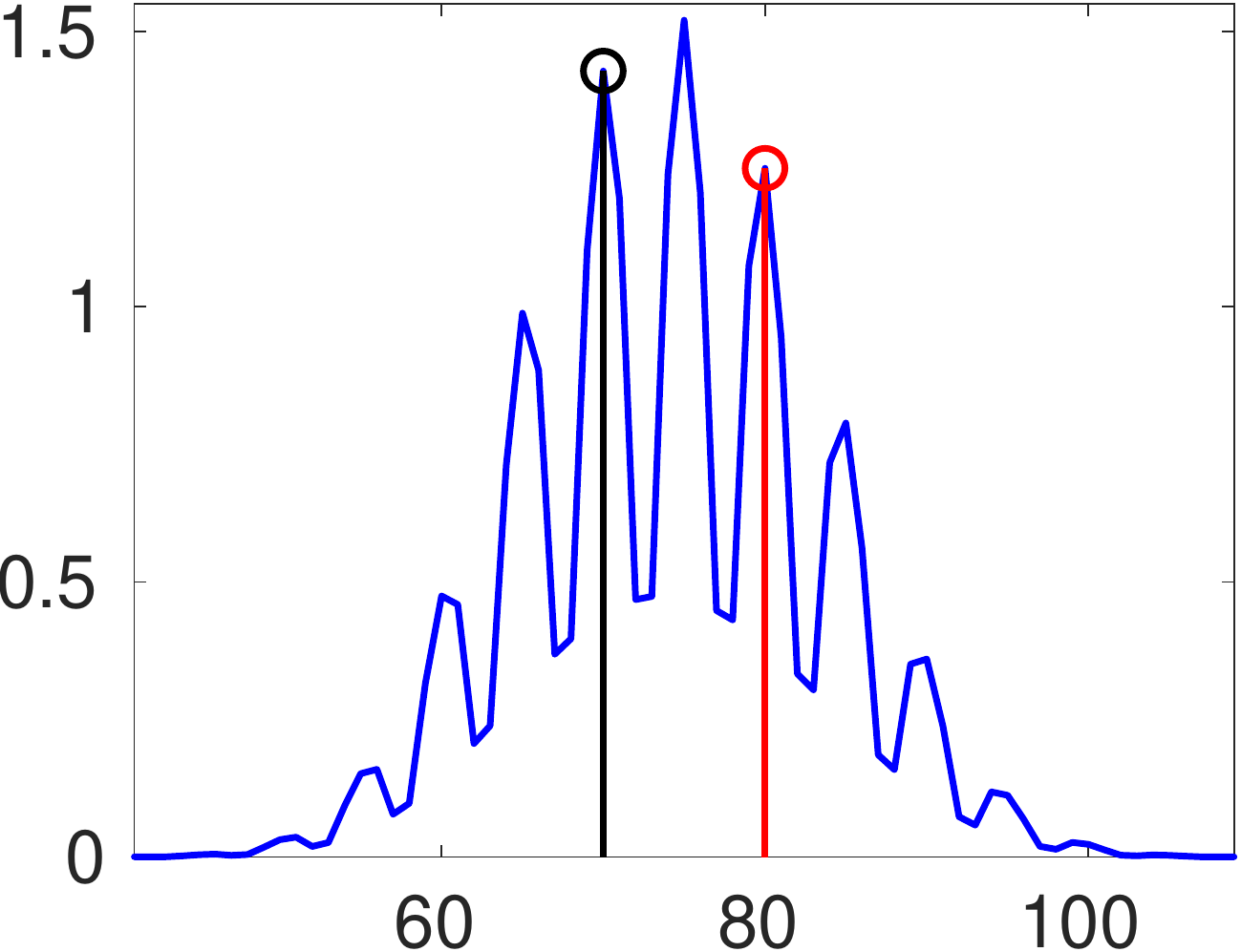}
		&	\includegraphics[width=\wth]{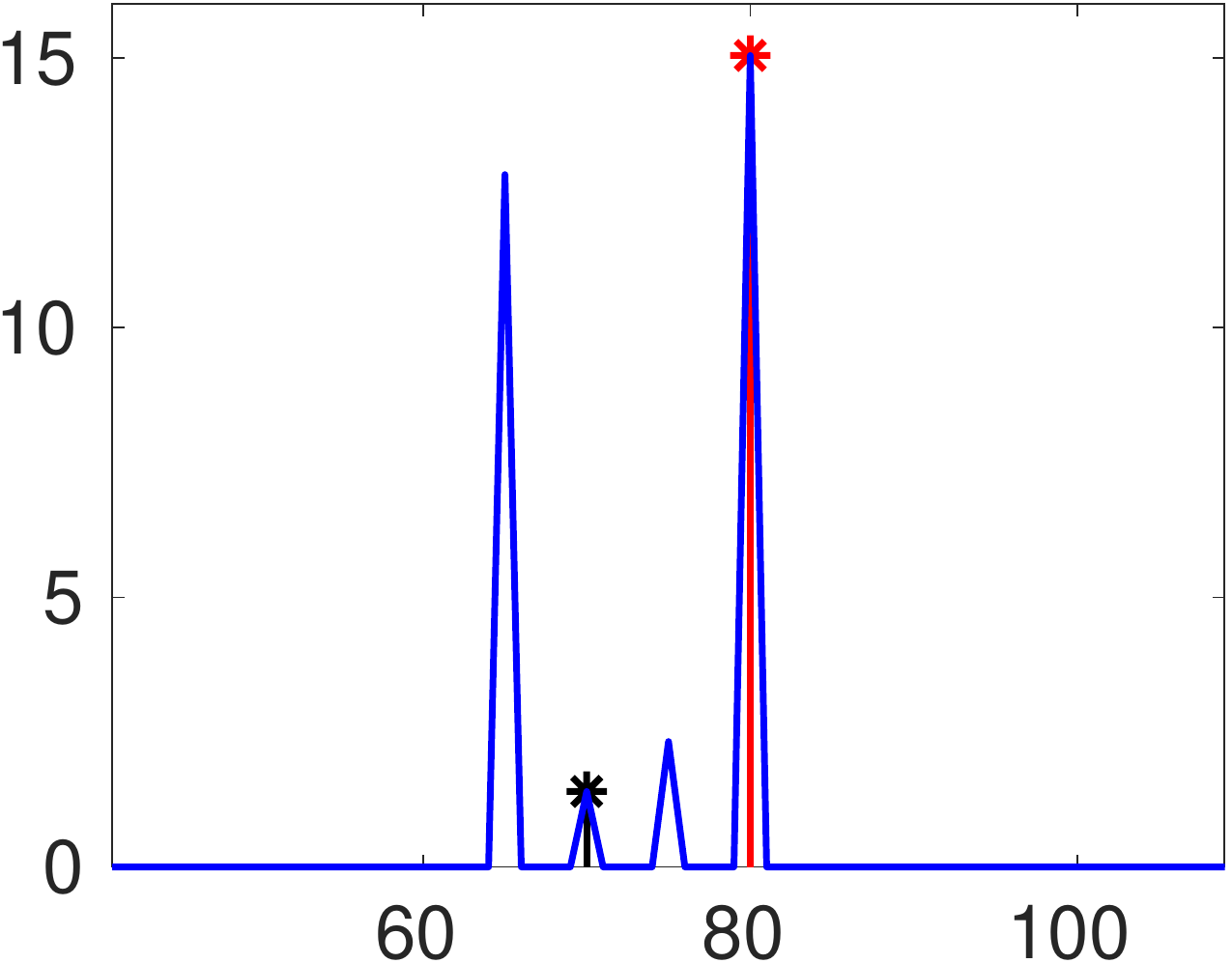}
		\\
		(b) & (d) & (f)\\
		\includegraphics[width=\wth]{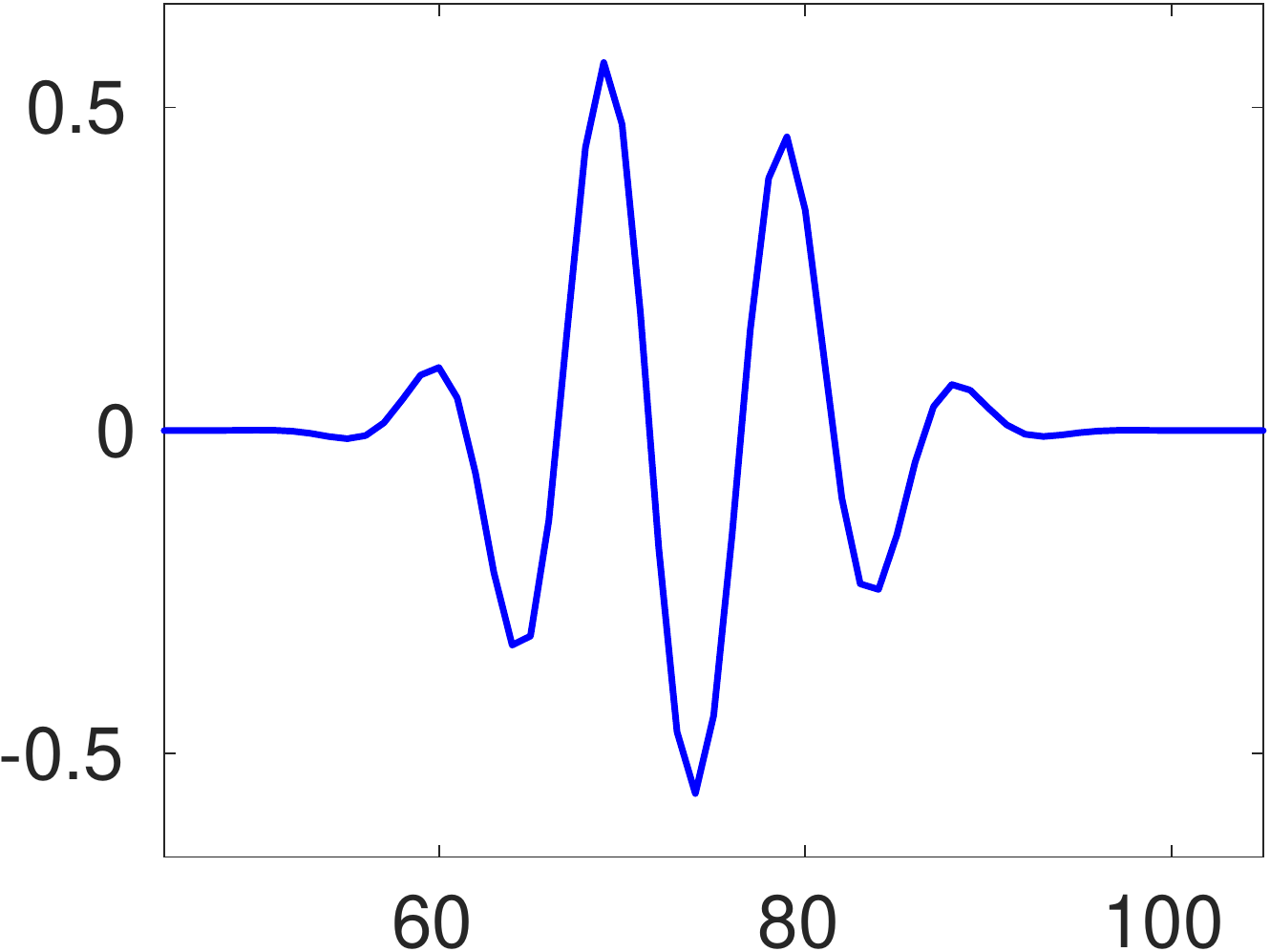}
		&	\includegraphics[width=\wth]{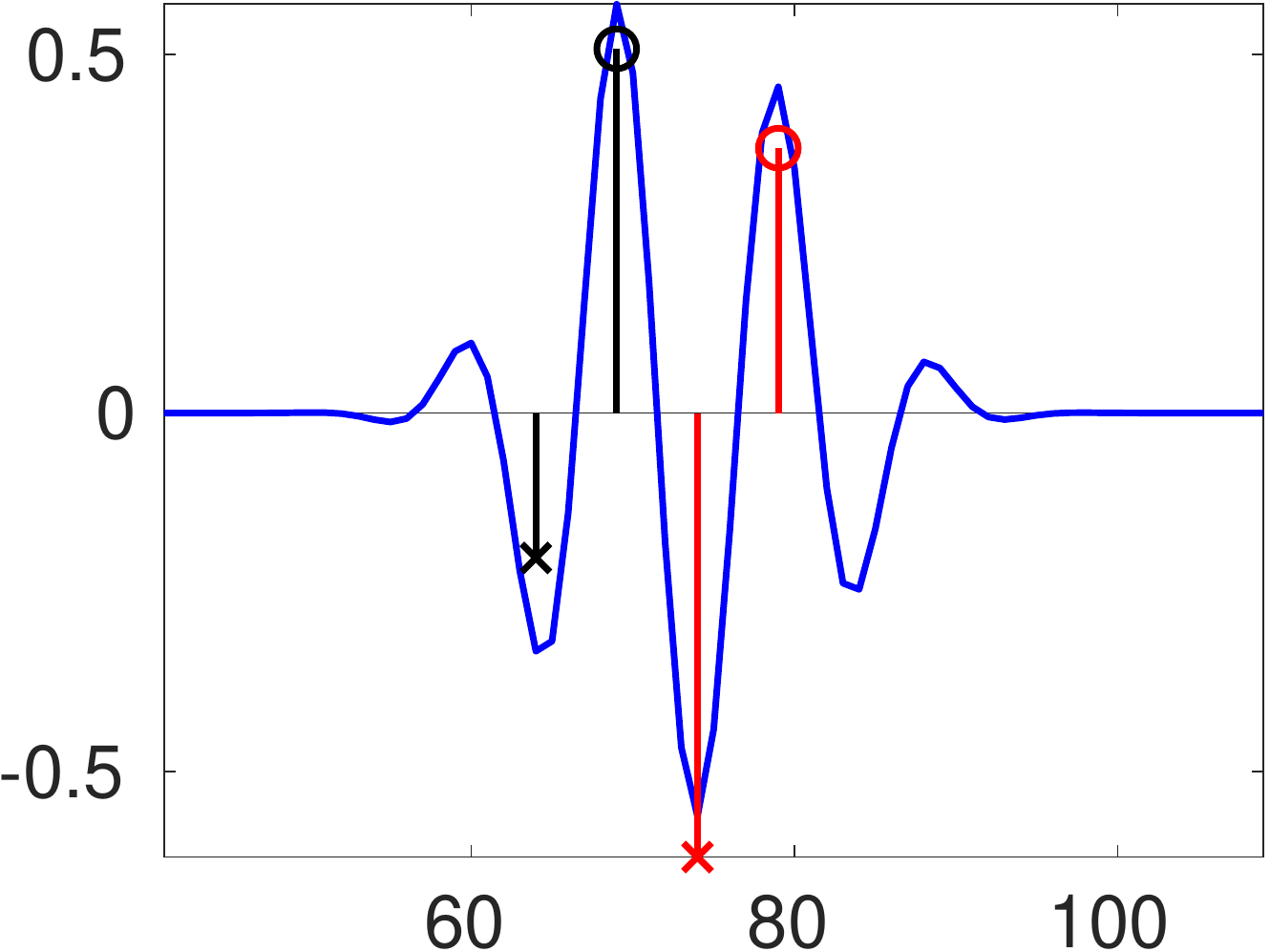}		
		&	\includegraphics[width=\wth]{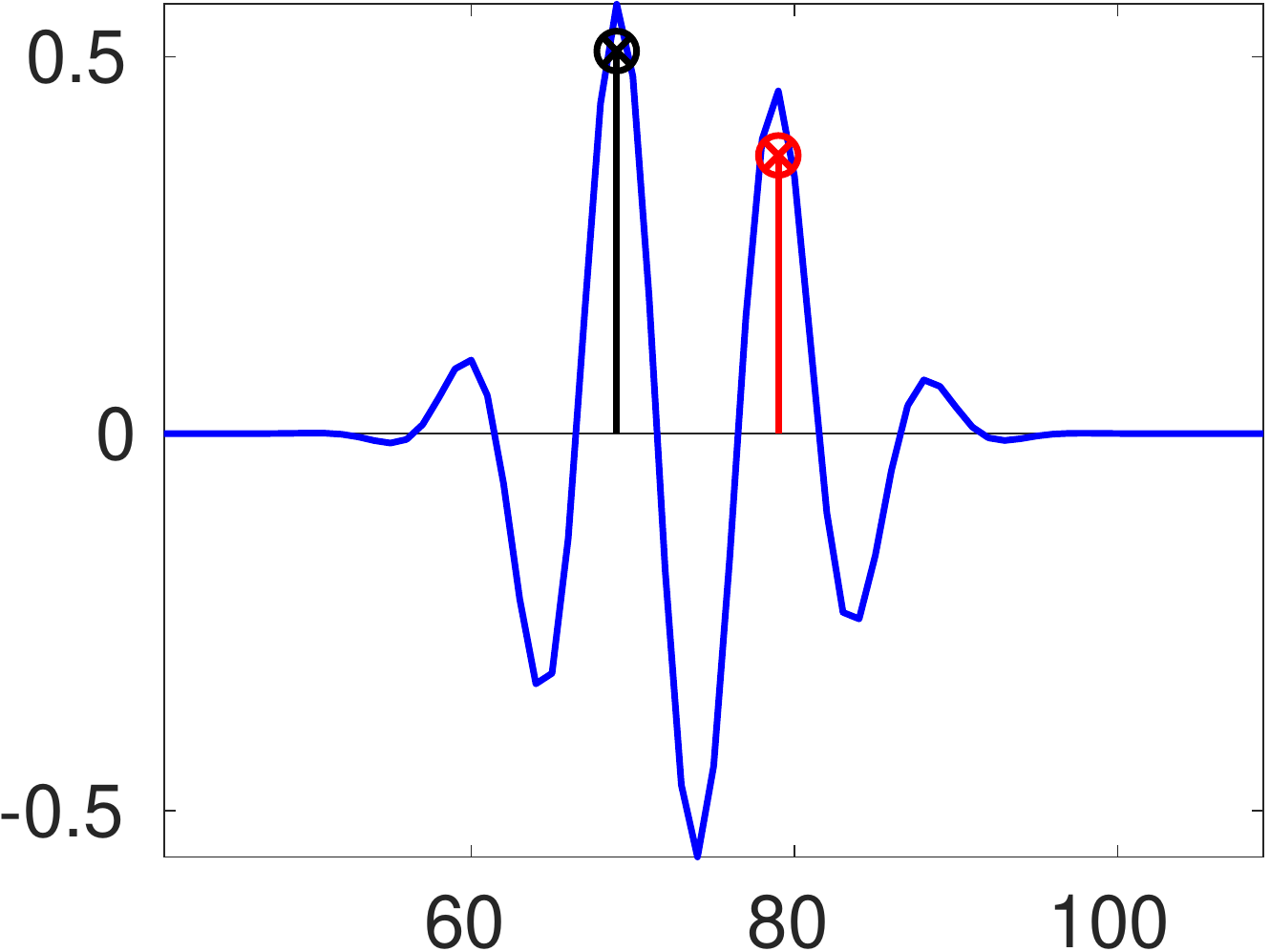}
	\end{tabular}
	\caption{a) Toy example where standard greedy selection fails. (a)~Contribution of two atoms and (b)~corresponding noise-free data. (c)~Scoring function for OMP and OLS at their first iteration and (d)~OLS solution after two iterations. (e)~Scoring function at the first iteration of the proposed SLS algorithm and (f)~SLS solution after two iterations. Circles (resp. stars) locate the true (resp. estimated) spikes.}
	% ({\color{magenta}\small $\triangle$}) and OMP ({\color{red}\bf$\circ$}). b) L1 ({\color{darkgreen}\small $\times$}) solution and true solution $\vx^\text{true}$ ({\color{black}\bf$\circ$}).}
	\label{fig::ex}
\end{figure}

\section{Single $\ell_1$ Selection}
\label{sec:SLS}

%Our selection rule exploits the solutions of $\ell_1$-norm-based problems. Recall the example in Figure\,\ref{fig::ex}. It is clear that a joint approach to the sparse estimation problem would be less sensitive to interferences between the atom contributions, by allowing non-zero weights to the true atom locations. 

We propose to select a new atom by considering the following $\ell_1$-norm optimization problem at each iteration:
\begin{equation}
{\min_{	\vx_{\Qup} \in \R^{{\nup}},
		\vx_{\Qbar} \in \mathbb{R}^{\barn} }} {\tfrac12 {\norm{\vy-\mH_{\Qup} \vx_{\Qup} - \mH_{\Qbar}\,\vx_{\Qbar}}}^2_2 + \lambda \norm{\vx_{\Qbar}}_1},
\label{eq:Plambda}
\end{equation}
Similarly to OLS, such a criterion allows the  re-estimation of the amplitudes of previously selected components $\vx_{\Qup}$ (whereas they are fixed in the residual for the selection rule of OMP) within the selection step. 
It also jointly estimates a sparse vector $\vx_{\Qbar}$ for the remaining ones, which is not restricted to a single non-zero component as in OLS.
Note that for a given $\vx_{\Qbar}$, the solution in $\vx_{\Qup}$ reads $\vx_{\Qup} = \mB^+   (\vy- \mH_{\Qbar} \vx_{\Qbar})$,
 therefore the problem in Eq.~\eqref{eq:Plambda} can be recast as an optimization problem in $\vx_{\Qbar}$ only:
% indeed, for a given $\vx_{\Qbar}$, the solution in $\vx_{\Qup}$ is explicit. %reads:
%$$	\vx_{\Qup} = {(\mB^T \mB)}^{-1} \mB^T (\vy- \mH_{\Qbar} \vx_{\Qbar}),$$
%Its expression can be inserted into the least-squares term in~\eqref{eq:Plambda}, so that the problem amounts to the following standard $\ell_1$-norm-penalized problem:
\begin{align}
%\begin{cases}
\widehat{\vx}_{\Qbar} = \arg\min_{\vx_{\Qbar}\in \mathbb{R}^{\barn}} \tfrac12 {\|\widetilde{\vy}- \widetilde{\mH}_{\Qbar} \,\vx_{\Qbar} \|}^2_2 + \lambda \|\vx_{\Qbar}\|_1%,\\
%\widehat{\vx}_{\Qup} = {(\mB^T \mB)}^{-1} \mB^T (\vy- \mH_{\Qbar} \widehat{\vx}_{\Qbar}),
\label{eq:pbl1_pen}
% \end{cases}
\end{align}
where $
\mP:= \mI_{\nup} - \mB \mB^+$, $\widetilde{\vy}:= \mP\vy$, and $\widetilde{\mH}_{\Qbar}:= \mP\mA$.% and  with $\mI_{\nup} $ the $\nup \times \nup$ identity matrix.

Then, the new component is selected as 
the component in $\widehat{\vx}_{\Qbar}$ with maximum amplitude value:
\begin{equation}
\score_\text{SLS}(j) = |\widehat{x}_j|, \; j \in \Qbar. \label{score_hom}
\end{equation}

In this paper, the solution to~\eqref{eq:pbl1_pen} is computed by the homotopy algorithm~\cite{Tibshirani96,donoho_fast_2008}, which iteratively solves~\eqref{eq:pbl1_pen} for decreasing values of parameter $\lambda$, starting from $\lambda_{\max}:= \max_{j \in \Qbar} |\widetilde{\va}_j^T \widetilde{\vy}|$, which is the minimum value of $\lambda$ above which the solution is identically zero. Here, we use a stopping rule that controls the number of nonzero components in $\widehat{\vx}_{\Qbar}$. More precisely, at a given iteration of the forward greedy procedure we impose that 
$$\norm{\widehat{\vx}_{\Qbar}}_0 =3 \,(K- \text{Card} (\Qup)),$$
where $K- \text{Card} (\Qup)$ represents the number of non-zero components that still must be determined. 
Such an empirical rule limits the computation time of the homotopy method, while ensuring that enough components are present in the model in order to overcome the interference issues explained in Section~\ref{sec:glouton}. 
%where $K-\nup$ is the number of components that still need to be included in the greedy search, and $c$ is a parameter whose tuning will be discussed in Section~\ref{sec:results}.

%Such an empirical rule follows the idea that a larger solution path is preferred for the first iterations of the greedy procedure, because more competition may exist between atoms, so the solution may be less stable as a function of $\lambda$.

%
The SLS selection rule is illustrated on the toy example of Section~\ref{sec:glouton} in Figure~\ref{fig::ex}. At the first iteration of the SLS algorithm,
the scoring function is maximum for a true component---it also shows non-zero, but smaller, value at the index erroneously selected by OMP and OLS, see Figure~\ref{fig::ex}\,(e). Then, after two iterations, the true support is correctly estimated, as shows Figure~\ref{fig::ex}\,(f). 
\section{Simulation results}
\label{sec:results}

We evaluate the performance of the SLS algorithm, compared to several well-known sparse estimation algorithms: OMP~\cite{Pati93},  OLS~\cite{Chen89}, $\ell_1$-norm regularization or BPDN (Basis Pursuit DeNoising), also computed here by the homopotopy algorithm~\cite{Tibshirani96}, SBR~\cite{soussen_bernoulli_2011}, Subspace Pursuit~\cite{Needell09}, accelerated Iterative Hard Thresholding (IHT)~\cite{blumensath_iterative_2008} and $A^\star$OMP~\cite{Karahanoglu12}.  
All algorithms are implemented in Matlab and are tuned such that all solutions have the true sparsity level $K$.

Algorithms are tested on difficult sparse deconvolution problems, with an up-sampled convolution model in order to achieve high-resolution spike locations~\cite{Carcreff13a}.
Problems are underdetermined with $\n =1,000$ and $\N = 350$. 
Columns of $\mH$ are then highly correlated, with mutual coherence $\max_{i\neq j} |\vh_i^T \vh_j| = 0.81$.  White Gaussian noise $\vect{\epsilon}$ is then added with $\text{SNR}_{\text{dB}} = 10\log\frac{{\Vert\mH\vx\Vert}^2}{{\Vert\vect{\epsilon}\Vert}^2} = 20$\,dB. %The $\mH$ columns are normalized.
\begin{figure}[h]
	\centerline{\includegraphics[width=8cm]{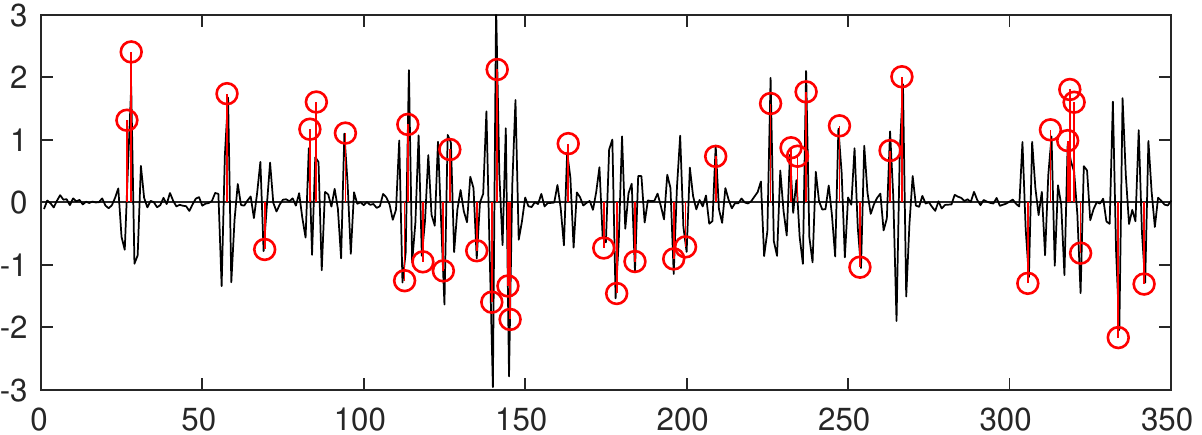}}
	\caption{Typical sparse deconvolution problem: sequence $\vx$ (\textcolor{red}{$\circ$}) and data~$\vy$~(-).}
	\label{fig:signaly}
\end{figure}
Results are averaged over 50 random realizations of the sparse sequence and of noise. A typical signal is shown in Figure~\ref{fig:signaly}.

\begin{figure}[h!]
	\centering
	\begin{tabular}{c@{\hspace*{0.5cm}}c}% @{\hspace*{0.1cm}}c @{\hspace*{0.1cm}}c @{\hspace*{0.1cm}}c}
	\includegraphics[width=4.1cm]{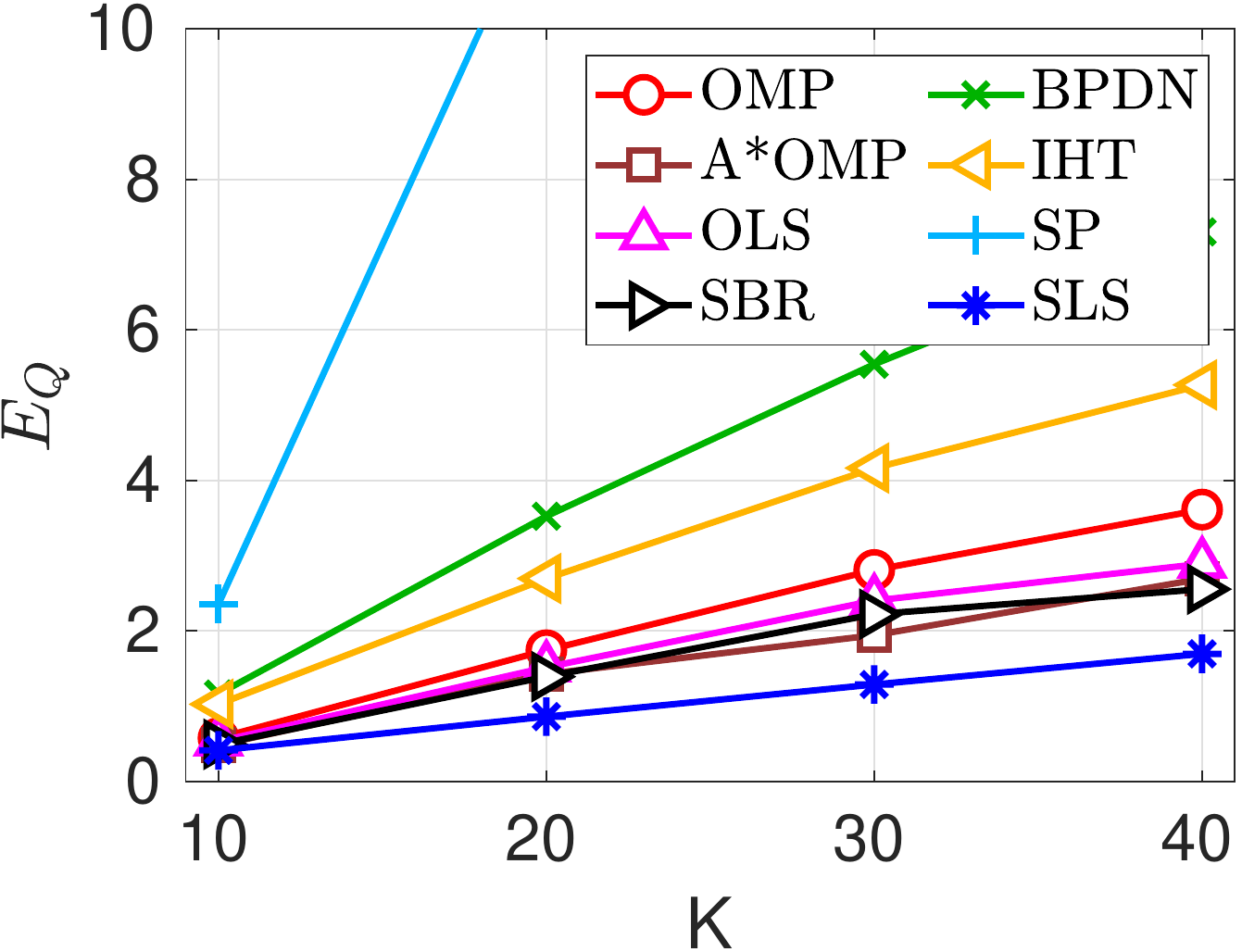}
	&\includegraphics[width=4.1cm]{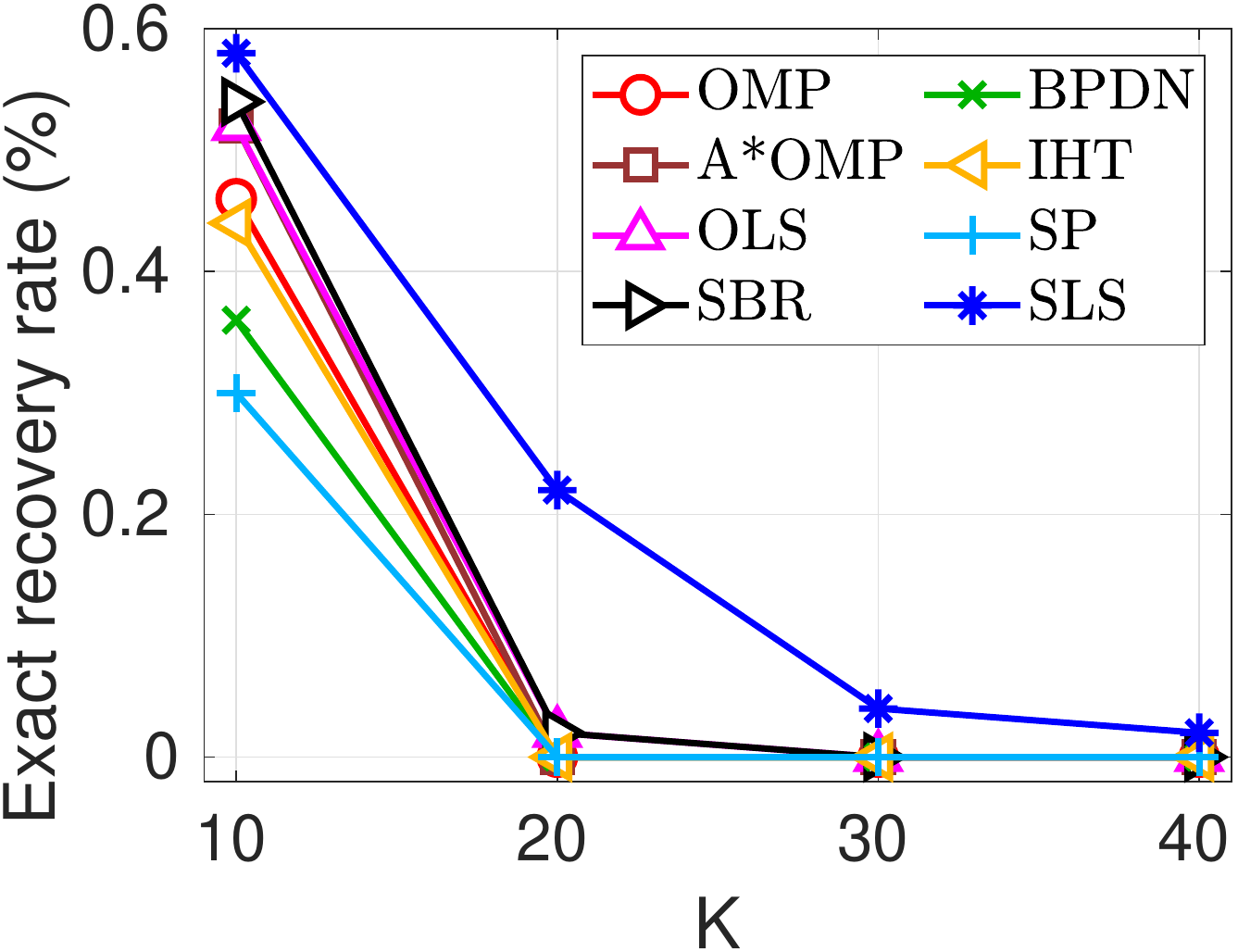}\\
	\multicolumn{2}{c}{\includegraphics[width=5cm]{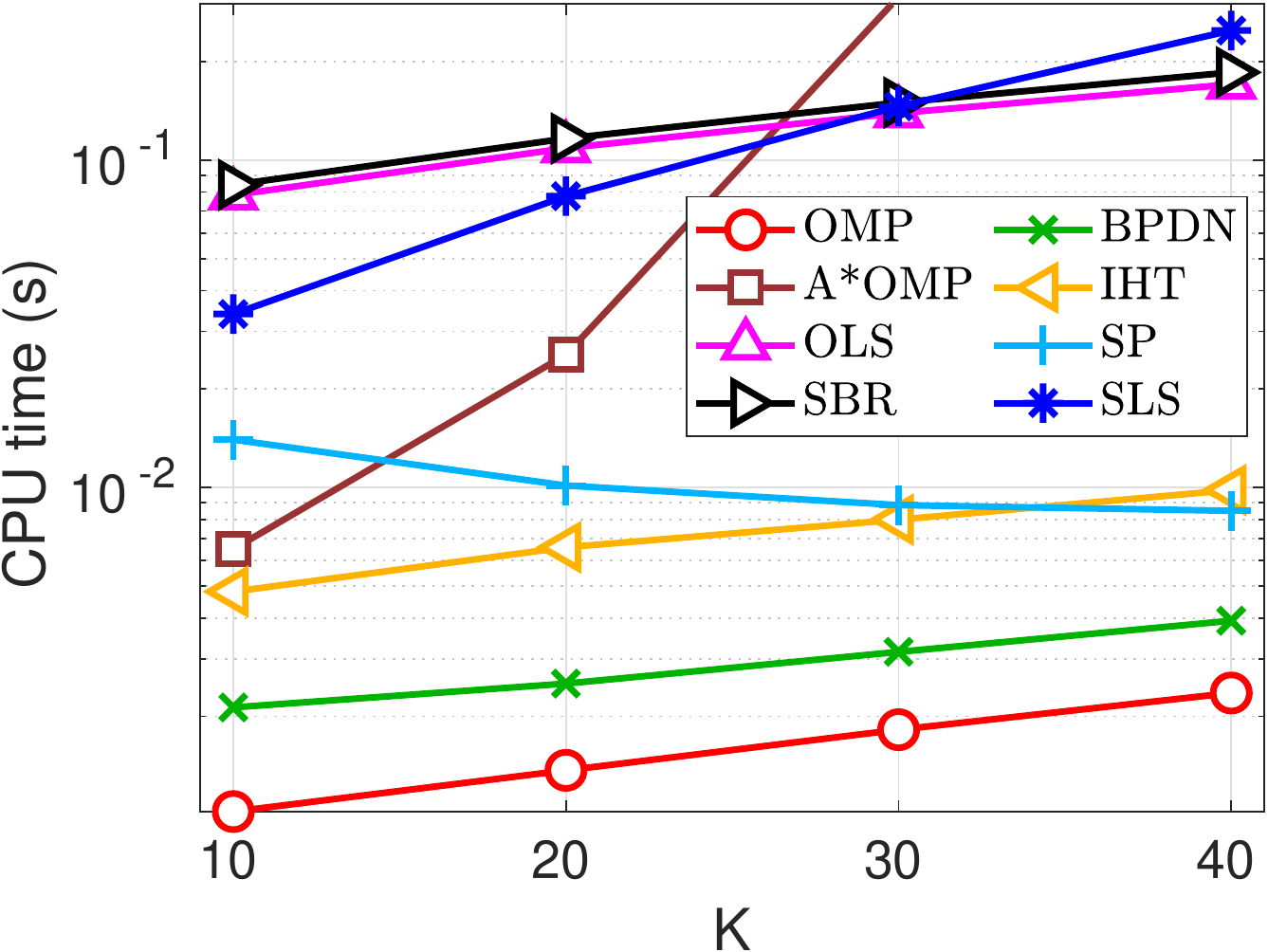}}\\
	\end{tabular}
	\caption{Quadratic error (top left), exact recovery rate (top right) and computation time as a function of $K$ for SLS and other sparse approximation methods.}
	\label{fig:res}
\end{figure}

Figure\,\ref{fig:res} shows the average quadratic error (left), the exact recovery rate (the fraction of solutions which have the correct support) and the average computing time for all algorithms as a function of $K$. 
SLS clearly achieves the best performance in terms of solution quality, and has a lower computation time than OLS and SBR up to $K=30$ and slightly higher for $K=40$. Note that other fast algorithms (OMP, BP, IHT and SP) are always much faster than SLS---but always give worse solutions.

% {\bf SLS vs OLS/SBR :} %SLS gives better solutions for different values of $K$ with 
% SLS has a better computation time up to $K=30$ and slightly higher for $K=40$.\\
%{\bf SLS vs $A^{\star}$OMP :} SLS %always gives better solutions for different values of K and it
%SLS has a slightly higher computation time up to $K=20$ but it becomes much faster from $K=30$.\\
%{\bf SLS vs OLS/SBR :} SLS %gives better solutions for different values of $K$ with a better computation time up to $K=30$ and slightly higher for $K=40$.

%\addRM{Sur la figure 4, il faudra quand même nuancer avec le coût de calcul. Mais dire en effet que l'algo fournit les meilleurs résultats (pour un coût de calcul inférieur jsuqu'à K= 30, et slightly supérieur ensuite).}

\bibliographystyle{IEEEbibJI}
\bibliography{bib_ICASSP}%strings,refs}

\begin{thebibliography}{10}

\bibitem{Taylor79}
H.~Taylor, S.~Banks, and F.~McCoy,
\newblock ``Deconvolution with the $l_1$ norm,''
\newblock {\em Geophysics}, vol. 44, no. 1, pp. 39--52, Jan. 1979.

\bibitem{Mendel83}
J.~Mendel,
\newblock {\em Optimal seismic deconvolution: {A}n estimation-based approach},
\newblock Monograph Series. Academic Press, 1983.

\bibitem{Zala92}
C.~A. Zala,
\newblock ``High-resolution inversion of ultrasonic traces,''
\newblock {\em IEEE Transactions on Ultrasonics, Ferroelectrics, and Frequency
  Control}, vol. 39, no. 4, pp. 458--463, July 1992.

\bibitem{Eldar12book}
Y.~Eldar and G.~Kutyniok,
\newblock {\em Compressed Sensing: Theory and Applications},
\newblock Cambridge University Press, 2012.

\bibitem{Miller02}
A.~Miller,
\newblock {\em Subset selection in regression},
\newblock Chapman and Hall/CRC, 2002.

\bibitem{natarajan_sparse_1995}
B.~K. Natarajan,
\newblock ``Sparse approximate solutions to linear systems,''
\newblock {\em SIAM Journal on Computing}, vol. 24, no. 2, pp. 227--234, 1995.

\bibitem{Tropp10}
J.~A. Tropp and S.~J. Wright,
\newblock ``{Computational methods for sparse solution of linear inverse
  problems},''
\newblock {\em Proceedings of the IEEE}, vol. 98, no. 6, pp. 948--958, June
  2010.

\bibitem{Pati93}
Y.~Pati, R.~Rezaiifar, and P.~S. Krishnaprasad,
\newblock ``Orthogonal matching pursuit: recursive function approximation with
  applications to wavelet decomposition,''
\newblock in {\em Asilomar Conference on Signals, Systems and Computers}, 1993,
  pp. 40--44 vol.1.

\bibitem{Chen89}
S.~Chen, S.~Billings, and W.~Luo,
\newblock ``Orthogonal least squares methods and their application to
  non-linear system identification,''
\newblock {\em International Journal of Control}, vol. 50, no. 5, pp.
  1873--1896, 1989.

\bibitem{soussen_bernoulli_2011}
C.~Soussen, J.~Idier, D.~Brie, and J.~Duan,
\newblock ``From {Bernoulli} {Gaussian} {Deconvolution} to {Sparse} {Signal}
  {Restoration},''
\newblock {\em IEEE Transactions on Signal Processing}, vol. 59, no. 10, pp.
  4572--4584, 2011.

\bibitem{Tibshirani96}
R.~Tibshirani,
\newblock ``Regression shrinkage and selection via the lasso,''
\newblock {\em Journal of the Royal Statistical Society, Series B}, vol. 58,
  pp. 267--288, 1996.

\bibitem{donoho_fast_2008}
D.~L. Donoho and Y.~Tsaig,
\newblock ``Fast {Solution} of $\ell_1$ {Norm} {Minimization} {Problems} {When}
  the {Solution} {May} {Be} {Sparse},''
\newblock {\em IEEE Transactions on Information Theory}, vol. 54, no. 11, pp.
  4789--4812, Nov. 2008.

\bibitem{Needell09}
D.~Needell and J.~Tropp,
\newblock ``Cosamp: Iterative signal recovery from incomplete and inaccurate
  samples,''
\newblock {\em Applied and Computational Harmonic Analysis}, vol. 26, no. 3,
  pp. 301 -- 321, 2009.

\bibitem{blumensath_iterative_2008}
T.~Blumensath and M.~E. Davies,
\newblock ``Iterative {Thresholding} for {Sparse} {Approximations},''
\newblock {\em Journal of Fourier Analysis and Applications}, vol. 14, no. 5,
  pp. 629--654, Dec. 2008.

\bibitem{Karahanoglu12}
N.~B. Karahanoglu and H.~Erdogan,
\newblock ``A* orthogonal matching pursuit: Best-first search for compressed
  sensing signal recovery,''
\newblock {\em Digital Signal Processing}, vol. 22, no. 4, pp. 555 -- 568,
  2012.

\bibitem{Carcreff13a}
E.~Carcreff, S.~Bourguignon, J.~Idier, and L.~Simon,
\newblock ``Resolution enhancement of ultrasonic signals by up-sampled sparse
  deconvolution,''
\newblock in {\em {{IEEE International Conference on Acoustic, Speech and
  Signal Processing}}}, Vancouver, Canada, May 2013.

\end{thebibliography}

\end{document}